\documentclass[11pt]{article}
\usepackage{latexsym,amsfonts,euscript}
\usepackage{xcolor}
\usepackage{amssymb}
\usepackage{stmaryrd}
\usepackage{mathrsfs}
\usepackage{leftidx}
\def\la{\langle}\def\ra{\rangle}

\def\pf{\noindent{\bf Proof\quad }}
\def\pfend{\hfill{$\Box$}}

\newtheorem{lem}{ \bf Lemma}[section]

\newtheorem{pro}[lem]{\bf Proposition}
\newtheorem{thm}[lem]{\bf Theorem}
\newtheorem{rem}[lem]{\bf Remark}

\newtheorem{cor}[lem]{\bf Corollary}

\makeatletter 
\@addtoreset{equation}{section}
\makeatother  

\title{On $p$-parts of  character degrees
\thanks{Project supported by
the  Nature Science Foundation of China (No. 11871011) and a grant from the Simons Foundation (No. 499532).   }}

\author{Guohua Qian\\
{\footnotesize\small  Dept. Mathematics, Changshu Institute of
Technology, Changshu, Jiangsu, 215500, China}\\
{\footnotesize\small E-mail: ghqian2000@163.com}\\\\
 Yong Yang\\
{\footnotesize\small  Dept. Mathematics, Texas State University, San Marcos, TX 78666, USA}\\
{\footnotesize\small E-mail: yang@txstate.edu}}

\begin{document}
\maketitle
\date{}

\vskip 1cm

\textbf{Abstract}\,\, Let $G$ be a finite group and $p$ be a prime.
In this paper,
we get the sharp bound for
$|G/O_p(G)|_p$ under the assumption that either $p^2 \nmid \chi(1)$ for all $\chi \in {\rm Irr}(G)$ or $p^2 \nmid \phi(1)$ for all $\phi \in {\rm IBr}_p(G)$. This would settle two conjectures raised by Lewis, Navarro, Tiep, and Tong-Viet in \cite{LNTV}.


\bigskip

\textbf{Keywords}\,\, $p$-parts of degree, character, Brauer character, Sylow $p$-subgroup.

\textbf{2010 MR Subject Classification}\,\,20C15 (primary),
20C20, 20D05 (secondary).

\vskip 5cm

\pagebreak

\section{Introduction}

In this paper,
$G$  always denotes  a finite group, and
${\rm Irr}(G)$ is the set of irreducible complex characters of $G$.
For a given prime $p$,
let $m_p$ be the $p$-part of an integer $m$, and let
$$e_p(G)=\max\{\log_p (\chi(1)_p)\,|\,\chi\in {\rm Irr}(G)\}.$$

\noindent
It is natural to ask how does this quantity $e_p(G)$ affect the structure of the Sylow $p$-subgroups of $G$. For example, the famous Ito-Michler Theorem \cite{Michler} asserts that $e_p(G)=0$ if and only  if $G$ has a normal abelian Sylow
$p$-subgroup. In particular, when this happens, we have
$|G/O_p(G)|_p=1$.

Let $b(P)$ denote the largest degree of an irreducible character of $P$. Moret\'o \cite[Conjecture 4]{Moret1} conjectured that $\log_p b(P)$ is bounded by some function of $e_p(G)$. The conjecture was first proved for solvable groups by Moret\'o and Wolf using an orbit theorem of solvable linear groups in \cite{MOWOLF}, and later settled by the authors for arbitrary finite groups in \cite{YQ}. In fact, we showed that $\log_p |G/O_p(G)|_p$ is bounded by a linear function of $e_p(G)$.

The research on the structure of groups under the condition $e_p(G)\leq 1$ was first studied in \cite{LNW}. In that paper, Lewis, Navarro, and Wolf showed that for solvable groups, if $e_p(G)\leq 1$, then $|G/O_p(G)|_p\leq p^2.$ For arbitrary finite groups, Lewis, Navarro, Tiep, and Tong-Viet proved in \cite{LNTV} that if $e_p(G)\leq 1$, then $|G/O_p(G)|_p\leq
p^4$. This result was improved to $|G/O_p(G)|_p\leq
p^3$ in \cite{Q18} by the first author.

The alternating group ${\rm Alt}_7$ shows that when
$e_2(G)=1$, the sharp bound for $|G/O_2(G)|_2$ is $2^3$. Also,
\cite[Example 3.3]{LNW} shows that  for a solvable group $G$ with
$e_p(G)=1$, the sharp bound for $|G/O_p(G)|_p$  is $p^2$. So it is natural to conjecture, as has been asked in \cite[Page 485]{LNTV}, \begin{center}
 \emph{If $e_p(G)=1$ for an odd prime $p$,
 whether the sharp bound for $|G/O_p(G)|_p$ is  $p^2$?
 }\end{center}
The following result gives a positive answer to the question along with some extra structural information.

\begin{thm}\label{t101} Let $G$ be a finite group
and $p$ be a prime.
If $e_p(G)\leq 1$, then $|G/O_p(G)|_p\leq p^2$, unless $p=2$ and $G\cong {\rm Alt}_7\times L$, where $L$ has a normal abelian Sylow $2$-subgroup, and in this case, it is clear that $|G/O_2(G)|_2=2^3$.


\end{thm}


Now we turn to consider
the similar question for
$p$-Brauer characters, which has also been considered in \cite{LNTV}.
As usual,
let ${\rm IBr}_p(G)$
be the set of irreducible $p$-Brauer characters of $G$.
Set
$$e'_p(G)=\max\{\log_p (\phi(1)_p)\,|\,\phi\in {\rm IBr}_p(G)\}.$$
A  Brauer character version of the Ito-Michler theorem only holds
for the given prime $p$, which asserts that $e'_p(G)=0$ if and only if $G$
has a normal Sylow $p$-subgroup (\cite[Theorem 5.5]{Michler}). For a
finite group $G$ with $e'_p(G)\leq 1$, it is proved in \cite[Theorem
1.5]{LNTV} that $$ |G/O_p(G)|_p\leq \left \{\begin{array}{lrrrrrrlr}
p^9, &{\rm if\,\, }&  p=2,\\
p^4, &{\rm if\,\, }&  p\geq 5 {\rm \,\,or,\,\, if}\,\,
p = 3 {\rm \,\, and\,\,}
{\rm Alt}_7 \,\,{\rm is\,\, not\,\, involved \,\,in}\,\,
G,\\
p^5, &{\rm if\,\,}& p = 3 {\rm \,\, and\,\,}
 {\rm Alt}_7 \,\,{\rm is\,\, involved \,\,in}\,\,
G.
\end{array}\right.$$
It is also conjectured \cite[Page 485]{LNTV} that the correct bounds
in the above inequalities are $p^7, p^2$ and $p^3$  respectively. We note that the bounds for the odd primes in the above inequalities ($p^4$ and $p^5$) have been slightly improved to $p^3$ and $p^4$ respectively in \cite{LY}. Generally speaking, ``$e'_p(G)\leq 1$'' does not imply ``$e_p(G)\leq
1$''.
Nevertheless, applying \cite[Theorem 1.1]{KM},
 which is the ``if part'' of the famous Brauer's height zero conjecture,
the authors in \cite{LNTV} have proved that if $G$ has an abelian
Sylow $p$-subgroup  with $O_p(G)=1$, then ``$e'_p(G)\leq 1$'' does
imply  ``$e_p(G)\leq 1$''. Note that if $G$ is  $p$-solvable with
$O_p(G) = 1$ and $e'_p(G)\leq 1$, \cite[Corollary 2.6]{W} shows that
$G$ has an elementary abelian Sylow $p$-subgroup. It follows that if
$G$ is $p$-solvable with $O_p(G)=1$, then $e'_p(G)\leq 1$ if and
only if $e_p(G)\leq 1$ (\cite[Corollary 1.4]{LNTV}). Using this
together with Theorem \ref{t101}, we obtain the upper bound for
$|G/O_p(G)|_p$ when $G$ is $p$-solvable with $e'_p(G)\leq 1$.

\begin{cor}\label{c102}
Let $G$ be a $p$-solvable group  for a prime $p$. If  $e'_p(G)\leq
1$, then $|G/O_p(G)|_p\leq p^2$.  \end{cor}

Assume that $G$ is non-$p$-solvable   with $e'_p(G)\leq 1$. For the
$p=2$ case, we  show that $G/O_2(G)$ is a direct product of $M_{22}$
and an odd order group (Proposition \ref{p304} below), and this  gives an
answer to the open questions proposed at the end of \cite{LNTV}.
For the case when $p\geq 3$, we show that $G/O_p(G)$ has an
elementary abelian Sylow $p$-subgroup (Proposition \ref{p305}),
whence $e_p(G/O_p(G))\leq 1$. Using these facts together with
Theorem \ref{t101}, we also get the upper bound for $|G/O_p(G)|_p$ along with some extra structural information.

\begin{thm}\label{t103}
Let $G$ be a finite group and $p$ be a prime.
If  $e'_p(G)\leq 1$,
then $|G/O_p(G)|_p\leq p^2$,
unless $p=2$ and $G/O_2(G)$ is a direct product of $M_{22}$ and an odd order group, in which case $|G/O_2(G)|_2=2^7$.
\end{thm}

We remark that
the bounds for $|G/O_p(G)|$ obtained in Theorems \ref{t101} and \ref{t103} are sharp.

We now fix some notation that will be used in the paper.

Let $G^{\sharp}$ be the set of nonidentity elements of $G$,
let ${\rm Irr}^{\sharp}(G)$ and ${\rm IBr}_p^{\sharp}(G)$ be
the set of nonprincipal irreducible characters and nonprincipal
irreducible $p$-Brauer characters, respectively, of $G$.

Let $E(p^n)$ be an elementary abelian $p$-group of order $p^n$. Let
$C_n$  denote a cyclic group  of order $n$.

Let ${\rm Sol}(G)$ be the largest solvable normal subgroup of $G$;
let ${\rm Soc}(G)$ be the socle of $G$, that is, the product of all
minimal normal subgroups of $G$.

For a given prime $p$, we denote by $G^{{\rm Sol}_p}$ the smallest normal subgroup of $G$
such that $G/G^{{\rm Sol}_p}$ is $p$-solvable,
and denote by ${\rm Sol}_p(G)$ the
largest normal $p$-solvable subgroup of $G$.

The paper is organized as follows. We will prove Theorem \ref{t101} in section 2 and Theorem \ref{t103} in section 3. 

\section{On $p$-parts of character degrees}

We will prove Theorem \ref{t101} at the end of this section. We begin with listing some known results
about finite groups $G$ with $e_p(G)\leq 1$.

\begin{lem} \label{l201}
Let $G$ be a finite group and
assume that $e_p(G)\leq 1$ for a prime  $p$. Then

{\rm (1)}  If $G$ is  solvable,
then $|G/O_p(G)|_p\leq p^2$;

{\rm (2)} If $p=2$ and $G$ is nonsolvable,
then $G$ is a direct product of ${\rm Alt}_7$ and
a finite group with  a normal abelian Sylow $2$-subgroup, in particular
$|G/O_2(G)|_2=2^3$;

{\rm (3)} If $p=3$ and ${\rm Alt}_7$ is involved in $G$, then
$O^{3'}(G)/O_3(G)\cong {\rm Alt}_7$, in particular $|G/O_3(G)|_3=3^2$;

{\rm (4)} If $p\geq 5$ or, if $p=3$ but ${\rm Alt}_7$ is not involved in $G$,
 then $|G/{\rm Sol}(G)|_p\leq p$.
\end{lem}

\pf  These are \cite[Theorem 1]{LNW}, \cite[Main Theorem]{L},
 \cite[Corollary 3.7]{LNTV} and  \cite[Lemma 3.3]{Q18}, respectively. \pfend

\bigskip

The results in the above lemma show that,
in order to prove Theorem \ref{t101},
we only need to consider the case
when $p\geq 3$ and $G$ is nonsolvable.

\begin{lem} \label{l202}  Suppose that
$W$ is a faithful and irreducible $G$-module for a solvable group
$G$. Assume that $G=O^{3'}(G)$, $|G|_3=3$ and $|W|\in \{ 2^2,2^4,
2^6,2^8\}$. Then $G\cong C_3$ and $W\cong E(2^2)$. \end{lem}

\pf This is a special case of \cite[Lemma 2.2]{LNW}. \pfend

\bigskip

Let $V$ be an irreducible $G$-module. Recall that  $V$  is called
imprimitive if $V$ can be written as $V=V_1\oplus \cdots \oplus V_n$
for $n\geq 2$  subspaces $V_i$ that are permuted transitively by
$G$, and that $V$ is primitive if $V$ is not imprimitive.

If $U\ltimes V$ is a Frobenius group with complement $U$ and kernel
$V$, then we say that $U$ acts fixed-point-freely on $V$.

\begin{lem}\label{l203}
Let $p,r$ be different primes with $p\geq 3$, let $G$ be a finite
group with $G=O^{p'}(G)$ that acts  faithfully and  completely
reducibly on an elementary abelian $r$-group $V$. Assume that
$|G|_p=p^c\geq p^2$, $C_p\cong O_p(G)$ acts fixed-point-freely on
$V$, and that every element of $V^{\sharp}$ is fixed by a unique
subgroup of order $p^{c-1}$ of $G$. Then $G$ acts primitively on
$V$.\end{lem}

\pf Let  $\mathfrak{X}$ be the set of subgroups $D$ of order
$p^{c-1}$ of $G$ such that $C_{V}(D)>1$. The hypothesis implies that
$$V = \bigcup_{D\in \mathfrak{X}} C_{V}(D),$$
\begin{equation}\label{eql20601} C_{V}(D_1)\cap C_{V}(D_2)
=\{1\}\end{equation}
whenever $D_1$ and $D_2$ are different members of $\mathfrak{X}$,
and that
\begin{equation}\label{eql20602}
O_p(G) \cap D=1
\end{equation}
for all $D\in \mathfrak{X}$. Assume that $G$ acts reducibly on $V$.
Write $V=W\times   U$, where $W$ and $U$ are nontrivial
$G$-invariant subgroups   of $V$. Let $w\in W^{\sharp}$, $u\in
U^{\sharp}$. Let  $D_*\in \mathfrak{X}$ be such that $D_*$
centralizes  $*$ where $*\in \{w,u, wu\}$. Since $W$ and $U$ are
$D_{wu}$-invariant, $D_{wu}$  must centralize both $w$ and $u$. Now
(\ref{eql20601}) yields
$$D_w=D_{wu}=D_u.$$ By the arbitrariness of $u$, $D_w$ centralizes
$U$. By the arbitrariness of $w$, $D_u$ centralizes $W$. This
implies that $D_w=D_u$  centralizes $V$ for all $w\in W^{\sharp}$
and all $u\in U^{\sharp}$. However $G$ acts faithfully on $V$, we
get a contradiction. Therefore, $G$ acts irreducibly on $V$.

Assume that $V$  is imprimitive. There exists subgroups  $V_1,
\ldots, V_d,\, d\geq 2$, of $V$ such that $V=V_1\times \cdots \times
V_d$ and $G$ acts transitively on $\{V_1, \ldots, V_d\}$. Let $K$ be
the kernel of the action. Note that $O_p(G)\cap K=1$ would imply
$C_V(O_p(G))>1$, this contradicts the assumption on $O_p(G)$.
Therefore $O_p(G)\leq K$. Clearly $G>K$, and thus $|G/K|_p\geq p$
because of $G=O^{p'}(G)$. Combining this together with
(\ref{eql20602}), we conclude that $K$ does not contain any subgroup
as a member of $\mathfrak{X}$.

Let $v_1\in V_1^{\sharp}$ and $D_1\in \mathfrak{X}$
be  such that $D_1$  centralizes $v_1$.
Since $D_1\not\leq K$, we may assume that $D_1$ does not normalizes  $V_2$.
Let $v_2\in V_2^{\sharp}$ and $D_2\in \mathfrak{X}$ so that
$D_2$ centralizes  $v_1v_2$. This implies that
$D_2$ acts on $\{v_1, v_2\}$.
Since $D_2$ is a $p$ group for an odd  prime $p$,
$D_2$ must centralize  both $v_1$ and $v_2$.
It follows by (\ref{eql20601})  that $D_1=D_2$.
Now $D_1$ centralizes $v_2$ and hence  normalizes
$V_2$, a contradiction.
Consequently $V$ is a primitive $G$-module, as wanted.\pfend

\bigskip
Let $V$ be an  irreducible $G$-module.
Recall that  $V$ is  quasi-primitive if $V_N$ is homogeneous
for all $N\unlhd G$. By Clifford's Theorem (\cite[Theorem 0.1]{MW}),
a primitive $G$-module is necessarily quasi-primitive.

Assume that a finite group $G$ acts on a set $\Omega$. For
$\alpha\in \Omega$,  $\alpha^G:=\{\alpha^g, g\in G\}$  is called a
$G$-orbit of $\alpha$. If $|\alpha^G|=|G|$, then the orbit
$\alpha^G$ is called regular; and if $|\alpha^G|_p=|G|_p$ for a
prime $p$, then the orbit $\alpha^G$ is called $p$-regular.

\begin{lem}\label{l204} Let $p, r$ be different primes with $p\geq 3$,
let  $G$ be  a solvable group with $G=O^{p'}(G)$, $E(p^2)\cong P\in
{\rm Syl}_p(G)$ and $O_p(G)\cong C_p$. Assume that $G$ acts
faithfully and quasi-primitively on an elementary $r$-group $V$. If
$F(G)$ is not cyclic and  $G$ has no $p$-regular orbit on $V$, then
$p=3$, $G\cong C_3\times {\rm SL}(2,3)$, and $V\cong E(7^2)$.
\end{lem}

\pf Clearly  $O_p(G)\cong C_p$ acts fixed-point-freely on $V$. Let
$\mathfrak{X}$ be the set of subgroups, different from $O_p(G)$, of
order $p$ of $G$. Since $G$ has no $p$-regular orbit on $V$, every
element of $V^{\sharp}$ is centralized by a  member of
$\mathfrak{X}$. This implies  that
\begin{equation}
\label{eql20403} V^{\sharp}=\bigcup_{D\in \mathfrak{X}} C_V(D)^{\sharp}.
\end{equation}
Since $V$ is a quasi-primitive $G$-module, every normal abelian
subgroup of $G$ is cyclic. Furthermore, since  $O^{p'}(G)=G$ and
$O_p(G)\cong C_p$, it is easily verified that every normal cyclic
subgroup of $G$  is also central in $F(G)$. Let $F = F(G)$ and let
$Z$ be the socle of cyclic group $Z(F)$. By \cite[Corollary
1.10]{MW}, there exist $E\unlhd G$, $U\unlhd G$ such that
\begin{equation}\label{eql20405}F=EU,\, E\cap U=Z, \, U=C_F(E)=Z(F),\,|F/U|=e^2;\end{equation}
$F/U$  is a completely reducible $G/F$-module and a faithful
$C_G(Z)/F$-module, possibly of mixed characteristic;  furthermore
$F/U$ is a direct product of irreducible $G/F$-modules of even
dimension.

Note that  ``$e=1$'' would imply $F(G):=F=U=Z(F)$ is  cyclic. Hence
$e>1$. Since $G$ has no  regular orbit on $V$, it follows by
\cite[Theorem 3.1]{Y15} that $e=2, 3, 4, 8, 9$ or $16$. Assume that
$p\geq 5$. The  proof of \cite[Theorem 3.2]{Y15} already show that
$G$ has a $p$-regular orbit on $V$, a contradiction. Therefore
$p=3$. Clearly  $$C_3\cong O_3(G)\leq Z\leq U,$$ and it follows that
$3\nmid e$. Consequently
$$p=3,\,e\in \{2,4,8,16\},\,6\mid |Z|,\,|G/F|_3=3.$$ Let
$K=C_G(F/U)$ and $F/U=V_1\times \cdots\times V_k$, where $V_1,
\ldots, V_k$ are irreducible $G/F$-modules. Since every $V_i$ has
even dimension,
 $|V_i|\in \{2^2,2^4, 2^6, 2^8\}$.
Obviously $|G/C_G(V_i)|_3=3$ because $O_3(G)\leq U\leq C_G(V_i)$.
Now  Lemma \ref{l202} yields $G/C_G(V_i)\cong C_3$. This also
implies that
$$G/K=G/(\bigcap_{i=1}^k C_G(V_i))
\leq G/C_G(V_1)\times \cdots \times G/C_G(V_k)$$
is a $3$-group.
Clearly $K<G$,
and since $F\leq K$ and $|G/F|_3=3$,
we conclude that
$$|G/K|=3, \,K/F=O_{3'}(G/F).$$
Since $F/U$ is a faithful $C_G(Z)/F$-module,
$K\cap C_G(Z)=F$.
By $G=O^{3'}(G)$ and the abelianlity of $G/C_G(Z)$,
$G/C_G(Z)$ is also a $3$-group. Consequently $C_G(Z)\geq K$.
Now  $F=K\cap C_G(Z)=K$ and thus
$$G/F\cong C_3.$$
For every  $D\in \mathfrak{X}$,
since $D\cap F=D\cap O_3(G)=1$, we have
\begin{equation}\label{eql20402} G=D\ltimes F.\end{equation}
Let $W$ be an irreducible constituent of $V_U$ where $U=Z(F)$. By
\cite[Theorem 2.2(7)]{Y11},
\begin{equation}\label{eql20406}
|V|=|W|^{eb}\end{equation} for a positive integer $b$.
Note that $6\mid |Z| \mid |U|$ and that $U$  acts fixed-point-freely on $W$.
It follows that
\begin{equation}\label{eql20407} |U|\mid (|W|-1),\, |W|\geq 7.\end{equation}
Applying (3) and (5) of \cite[Lemma 2.4]{Y10},
we conclude that
\begin{equation}\label{eql20408}
|C_V(D)|\leq |W|^{\frac{1}{2}eb}
\end{equation} for all $D\in \mathfrak{X}$.
Now  (\ref{eql20403}),(\ref{eql20406}) and (\ref{eql20408}) yield
\begin{equation}\label{eql20409}
|W|^{eb}-1=|V|-1\leq |\mathfrak{X}|(|W|^{\frac{1}{2}eb}-1).
\end{equation}
Clearly $|N_G(P|\geq |C_G(P)|\geq |U|_2\cdot|P|\geq 18$. Consequently
$$|G: N_G(P)|\leq \frac{|G:F||F|}{18}=\frac{1}{6}|F|=\frac{1}{6}|F:U||U|\leq \frac{1}{6}e^2(|W|-1).$$
Observe that $P\cong E(3^2)$ contains exactly
$3=\frac{3^2-1}{3-1}-1$ members of $\mathfrak{X}$, and it follows that
\begin{equation}\label{eql20410} |\mathfrak{X}|\leq 3|G: N_G(P)|
\leq \frac{1}{2}e^2(|W|-1).\end{equation}
By (\ref{eql20409}) and (\ref{eql20410}), we conclude that
$
|W|^{eb}-1\leq  \frac{1}{2}e^2(|W|-1)(|W|^{\frac{1}{2}eb}-1)$, that is,
\begin{equation}\label{eql20411}
|W|^{\frac{1}{2}eb}+1\leq \frac{1}{2}e^2(|W|-1),
\end{equation}
in particular,
$$
|W|^{\frac{1}{2}eb-1}< \frac{1}{2}e^2.$$
Since  $|W|\geq  7$ by (\ref{eql20407}) and  $e\in \{2,4,8,16\}$,
the above inequality  implies  $e\leq 4$ and $b=1$.
Note that if $e=4$ and $b=1$,
then (\ref{eql20411}) yields a contradiction.
Consequently,  $e=2$ and $b=1$.

Suppose that $U\not= Z(G)$ and let $D\in \mathfrak{X}$. Since $D$
acts nontrivially on the cyclic group $U$ by (\ref{eql20402}), there
exists a $3'$-subgroup $Y$ of $U$ of prime order such that $D$ acts
fixed-point-freely on $Y$. Considering the action of $D\ltimes Y$ on
$V$ and observing that $Y$ acts fixed-point-freely on $V$, we
conclude from  \cite[Lemma 0.34]{MW} that $\dim
C_V(D)=\frac{1}{3}\dim V$ for all $D\in \mathfrak{X}$. By
(\ref{eql20403}), (\ref{eql20406}) and (\ref{eql20410}), we get that
$$6\mid \dim V$$ and that
$$|V|-1\leq |\mathfrak{X}|(|V|^{1/3}-1)
\leq 2(|V|^{1/2}-1)(|V|^{1/3}-1)< 2(|V|^{5/6}-1),$$ which is clearly
impossible. Consequently $U=Z(G)$.

By (\ref{eql20405}), we have  $F/E\leq Z(G/E)$. Since $E\geq Z\geq
O_3(G)$, $F/E$ is a $3'$-group , and it follows that  $F=E$ because
$G=O^{3'}(G)$. Now
$$Z(G)=U=U\cap F=U\cap E=Z$$ via (\ref{eql20405}).
Recall that $F/U\cong E(2^2)$ because $e=2$. Obviously $G=O^{3'}(G)$
also implies $Z=Z(G)\cong C_6$.
 Now it is easy to see that
$$O_2(F)\cong Q_8,\,\,G=O_3(G)\times H, \,\,
{\rm where\,\,} H=D\ltimes O_2(F)\cong {\rm SL}(2,3)$$
for some   $ D\in \mathfrak{X}$, and that
$$|\mathfrak{X}|=3|G:N_G(P)|=12.$$
By (\ref{eql20409}), we get that $|W|^2-1\leq 12(|W|-1)$.
Hence   $|W|\leq 11$.
Since  $6=|U|$ divides $|W|-1$,
we have  $|W|=7$.
Recall that $|V|=|W|^{eb}$ where $e=2$ and $b=1$.
Consequently  $V\cong E(7^2)$, and we are done.   \pfend

\bigskip

In order to prove Theorem \ref{t206}, Propositions \ref{p304} and \ref{p305},
we also need the following lemma.

\begin{lem}\label{l205}
Let $P$ be a Sylow $p$-subgroup of a finite group $G$,
and let $Z\leq P\cap Z(G)$.
Assume that  $P$ is either abelian or split over $Z$.
Then $Z\cap G'=1$.
\end{lem}

\pf By the hypothesis, every $\lambda\in {\rm Irr}(Z)$ is extendible
to $P$. It follows by \cite[Theorem 6.26]{I} that $\lambda$ extends
to $\lambda_0\in {\rm Irr}(G)$. Observe that $\lambda_0$ is linear
and $G'\leq \ker\lambda_0$. We have  $\ker\lambda= Z\cap
\ker\lambda_0\geq Z\cap G'$. Therefore  $Z\cap G'\leq
\bigcap_{\lambda\in {\rm Irr}(Z)}\ker\lambda=1$. \pfend

\bigskip

\begin{thm}\label{t206} Let $H$ be a finite group and $p$ be an odd prime.
Assume that $e_p(H)\leq 1$. Then $|H/O_p(H)|_p\leq p^2$.\end{thm}

\pf Assume the result is not true.
Let $H$ be a counterexample with minimal order.
Note that the hypothesis $e_p(H)\leq 1$ is inherited
by all quotient groups and normal subgroups of $H$.
We will deduce  a contradiction via several steps.

\bigskip

(1) $O_p(H)=\Phi(H)=1$, $H=O^{p'}(H)$, $|H/{\rm Sol}(H)|_p=p$,
$|{\rm Sol}(H)|_p=p^2$ and $|H|_p=p^3$.

Suppose that $O_p(H)>1$.
Since $e_p(H/O_p(H))\leq 1$ and $H$ is a minimal counterexample,
we have $|(H/O_p(H))/O_p(H/O_p(H))|_p\leq p^2$,
that is, $|H/O_p(H)|_p\leq p^2$, a contradiction.
Therefore $O_p(H)=1$.

Suppose that $\Phi(H)>1$. Note that $\Phi(H)\leq O_{p'}(H)$ and
$O_p(H/\Phi(H))=1$ because $O_p(H)=1$. Since $e_p(H/\Phi(H))\leq 1$
and $H$ is a minimal counterexample, we conclude that
$|(H/\Phi(H))/O_p(H/\Phi(H))|_p\leq p^2$. This implies
$|H/O_p(H)|_p\leq p^2$, a contradiction. Therefore $\Phi(H)=1$.

Since $e_p(O^{p'}(H))\leq 1$, we also have  $H=O^{p'}(H)$.

By  Lemma \ref{l201}(3), we may assume that if $p=3$,
then ${\rm Alt}_7$  is not involved in $H$.
By (1) and (4) of Lemma \ref{l201},
we have that
$|H/{\rm Sol}(H)|_p\leq p$ and $|{\rm Sol}(H)|_p\leq p^2$.
Since $H$ is a counterexample,
we get that
$|H/{\rm Sol}(H)|_p=p$, $|{\rm Sol}(H)|_p=p^2$ and $|H|_p=p^3$.

\bigskip

(2) $H=G\ltimes V$, where $G$ is a maximal subgroup of $H$ with
$G=O^{p'}(G)$, $O_p(G)=\la a\ra\cong C_p$, $V\cong E(r^n)$ is a
unique minimal normal subgroup of $H$ where  $r\not=p$ is a prime,
furthermore $G$ acts faithfully on $V$.

Let $N$ be a minimal normal subgroup of $H$.
Assume that $p$ divides $|N|$.
Since $O_p(H)=1$,
we have $N\cap {\rm Sol}(H)=1$.
Observe that neither $N$ nor ${\rm Sol}(H)$
contains a normal abelian Sylow $p$-subgroup.
By the Ito-Michler Theorem,
there exist
$\theta_1\in {\rm Irr}(N)$ and  $\theta_2\in {\rm Irr}({\rm Sol}(H))$
such that $p$ divides $\theta_1(1)$ and $\theta_2(1)$.
Since $\theta_1\theta_2\in {\rm Irr}(N\times {\rm Sol}(H))$,
we get  that
$e_p(H)\geq e_p(N\times {\rm Sol}(H))\geq \log_p ((\theta_1\theta_2)(1))_p\geq 2$,
a contradiction.  Hence $N$ is necessarily a $p'$-group.
Note that since  $H$ is a minimal counterexample,
we also conclude that
$O_p(H/N)>1$ for every minimal normal subgroup $N$ of $G$.

Suppose that $H$ admits different minimal normal subgroups,
say $N_1, N_2$.
Let $D_i\unlhd H$ be such that $D_i/N_i=O_p(H/N_i)$.
Since $p\mid |D_i|$  and $O_p(D_i)\leq O_p(H)=1$ ,
there exists $\theta_i\in {\rm Irr}(D_i)$ such that $p\mid \theta_i(1)$
for every $i\in \{1,2\}$.
Note that $D_1\cap D_2=1$ and $D_1D_2=D_1\times D_2\unlhd H$.
It follows that $\theta_1\theta_2\in {\rm Irr}(D_1\times D_2)$ and that
$e_p(H)\geq \log_p ((\theta_1\theta_2)(1))_p\geq 2$,
a contradiction.
Hence $H$ admits a unique minimal normal subgroup, say $V$.

Since ${\rm Sol}(H)$ is nontrivial by (1), the uniqueness of $V$
implies $V\leq {\rm Sol}(H)$, and hence $V\cong E(r^n)$ for some
prime power $r^n$, where $r\not=p$. Since $\Phi(H)=1$, there exists a
maximal subgroup $G$ of $H$ such that $H=G\ltimes V$. Since
$H=O^{p'}(H)$, we have $G=O^{p'}(G)$.

Observe that $C_G(V)$ is normal in $H$.
The uniqueness of $V$ implies $C_G(V)=1$,
that is, $G$ acts faithfully on $V$.

Observe that $1<O_p(H/V)\cong O_p(G)\leq {\rm Sol}(H)$ and $|{\rm
Sol}(H)|_p=p^2$. It follows that  $O_p(G)$ has order $p$ or $p^2$.
In particular $O_p(G)$ is abelian. Obviously $O_p(G)$ acts
faithfully and coprimely on $V$. By \cite[Lemma 18.1]{MW}, there
exists $\theta\in {\rm Irr}(O_p(G)\ltimes V)$ of degree $|O_p(G)|$.
Since $e_p(O_p(G)\ltimes V))\leq e_p(H)\leq 1$, we have $O_p(G)=\la
a\ra \cong C_p$.

\bigskip

{\rm (3)} Let $P\in {\rm Syl}_p(G)$,
$J=O^{p'}({\rm Sol}(G))$ and $P_J=P\cap J$.
Let $J_0\in \{J, G\}$ and let $|J_0|_p=p^c$. Then

\quad (3a)
For every $\lambda\in {\rm Irr}^{\sharp}(V)$, we have that
${\rm I}_{J_0}(\lambda)\cap \la a\ra=1$ and that
${\rm I}_{J_0}(\lambda)$ admits a
normal abelian Sylow $p$-subgroup of order $p^{c-1}$;

\quad (3b) $P\cong E(p^3)$, $P_J\cong E(p^2)$;

\quad (3c) $J_0$ acts faithfully and  primitively on $V$ and ${\rm
Irr}(V)$.

Note that $|P|=p^3$ and $|P_J|=p^2$ by (1).

Clearly $\la a\ra=O_p(G)\leq J_0$, $J_0\unlhd G$ and thus
$e_p(J_0\ltimes V)\leq 1$. Let $\lambda\in {\rm Irr}^{\sharp}(V)$
and let $T={\rm I}_{J_0}(\lambda)$. Since $G$ acts faithfully and
irreducibly  on $V$, $\la a\ra$  acts  fixed-point-freely on $V$.
This implies that $T\cap \la a\ra=1$ and hence  $|T|_p\leq p^{c-1}$.
Assume that $|T|_p\leq p^{c-2}$. By \cite[Theorem 6.11]{I}, there
exists $\chi\in {\rm Irr}(J_0\ltimes V|\lambda)$ of degree divisible
by $p^2$, a contradiction. Assume that  $|T|_p=p^{c-1}$ and that  a
Sylow $p$-subgroup of $T$ is  nonnormal or nonabelian. Note that
$\lambda$ is extendible to $T\ltimes V$ and that $T$ has an
irreducible character of degree divisible by $p$. By \cite[Corollary
6.17]{I}, there exists  $\theta\in {\rm Irr}(T\ltimes V|\lambda)$ of
degree divisible by $p$; and by \cite[Theorem 6.11]{I}, there exists
$\chi\in {\rm Irr}(J_0\ltimes V|\theta)$ of degree divisible by
$p^2$, a contradiction. Consequently, $T$ admits a normal abelian
Sylow $p$-subgroup of order $p^{c-1}$.

Assume that $P/\la a\ra$ is cyclic. Since $O_p(G/\la a\ra)=1$, it is
well-known  \cite{Z} that $G/\la a\ra $ admits an irreducible
character of degree divisible by $|P/\la a\ra|=p^2$,  a
contradiction. Hence $P/\la a\ra\cong E(p^2)$. Let $\lambda\in {\rm
Irr}^{\sharp}(V)$ and $P_1\in {\rm Syl}_p({\rm I}_G(\lambda))$. Without
loss of generality, we may assume $P_1\leq P$. Note that  $\la
a\ra\cap P_1=1$ by (3a), $P_1\cong P/\la a\ra$ and $\la a\ra\unlhd
P$. We get that $P=\la a\ra \times P_1\cong E(p^3)$ and $P_J\cong
E(p^2)$.

Clearly $O^{p'}(J_0)=J_0$ and $J_0$ acts faithfully and completely
irreducibly on $V$ and ${\rm Irr}(V)$. It follows by (3a) and Lemma
\ref{l203} that ${\rm Irr}(V)$ is a faithful  primitive
$J_0$-module, and so is $V$. Note that we only use the primitivity
of $J$-module ${\rm Irr}(V)$ in the sequel.

\bigskip

(4) $F(J)$ is  cyclic.

Observe that   ${\rm Irr}(V)$  is a faithful  primitive $J$-module by (3),
and that $J$ has no $p$-regular orbit on ${\rm Irr}(V)$
 because  $e_p(J\ltimes V)\leq 1$ and $|P_J|=p^2$.
Assume that $F(J)$ is not cyclic. Applying Lemma \ref{l204},
we get that
$$p=3, J\cong C_3\times {\rm SL}(2,3),
V\cong {\rm Irr}(V)\cong E(7^2).$$ Since $V$ is the unique minimal
normal subgroup of $H$, we have $C_G(V)=1$ and $G\leq {\rm
GL}(2,7)$. Investigating the index of $J$ in ${\rm GL}(2,7)$, we
conclude that $G/J$ is a $\{2,7\}$-group. By Burnside's $p^aq^b$
Theorem, $G/J$ is solvable, and so is $H$, a contradiction.
Consequently $F(J)$ is cyclic.

\bigskip

(5) There exist elements $b, c\in P^{\sharp}$
and normal subgroups $B, C$ of $G$ such that

\quad (5a)  $P=\la a\ra \times \la b \ra \times \la c\ra$,
$G=\la a \ra \times (BC)$, where $\la a\ra=O_p(G)$;

\quad (5b) $\la a\ra \times B=J$,
$B=\la b\ra \ltimes O_{p'}(J)$
with $B=O^{p'}(B)$;

\quad (5c) $C$ is nonsolvable with $\la c\ra\in {\rm Syl}_p(C)$
and $C=O^{p'}(C)$,
$B\cap C=O_{p'}(B)=O_{p'}(J)$.

Since $C_p\cong \la a\ra< P_J$  and $P_J\cong E(p^2)$ by (3b), we
may write $P_J=\la a\ra \times \la b\ra$, where $o(a)=o(b)=p$. Since
$P_J\in {\rm Syl}_p(J)$  is  abelian, the  $p$-solvable subgroup $J$
with $J=O^{p'}(J)$ has a normal $p$-complement, that is,
$$J=P_J\ltimes O_{p'}(J).$$ Let  $M\lhd G$ be maximal such that
$P_J\cap M=1$ and $O_{p'}(J)\leq M$. We claim that $G=P_J\ltimes M$.
Observe that
$$E(p^2)\cong P_JM/M=P_JO_{p'}(J)M/M=JM/M\unlhd G/M.$$ Since $G/M$ admits an abelian Sylow $p$-subgroup,
$G/C_{G}(P_JM/M)$ is a $p'$-group. Since $G=O^{p'}(G)$, $P_JM/M$ lies
in the center of $G/M$. By Lemma \ref{l205}, we have $(P_JM/M)\cap
(G/M)'=1$, that is, $P_JM\cap G'M=M$. This also implies that
$$P_J\cap G'\leq P_J\cap G'M=P_J\cap (P_JM\cap G'M)=P_J\cap M=1.$$
By the maximality of $M$,
we have $G'M=M$ and thus $G'\leq M$.
Note that $G/M$ is an elementary abelian $p$-group
because $G=O^{p'}(G)$ and $E(p^3)\cong P\in {\rm Syl}_p(G)$.
Consequently $P_JM/M$ is split over $G/M$.
Let $U/M$ be a complement of $P_JM/M$ in $G/M$.
We have that
$$P_J\cap U=P_J\cap (P_JM\cap U)=P_J\cap M=1$$  and that
$G=P_JMU=P_JU=P_J\ltimes U$. Obviously $U=M$ and the claim follows.

Let $B=\la b\ra \ltimes O_{p'}(J)$. Clearly $J=P_J\ltimes O_{p'}(J)
=(\la a\ra \times \la b))\ltimes O_{p'}(J)=\la a\ra \ltimes B=\la
a\ra\times B$. Since $J=O^{p'}(J)$, we have $B=O^{p'}(B)$. Note that
$\la b\ra M \unlhd G$ because $G=P_J\ltimes M$ by the claim, and
that $$(P_J\cap \la b\ra M) O_{p'}(J)= P_JO_{p'}(J)\cap \la b\ra M
=J\cap \la b\ra M $$ where the first equality follows from
$O_{p'}(J)\leq M\leq \la b\ra M$. This implies that
$$B=\la b\ra \ltimes O_{p'}(J)=(P_J\cap \la b\ra M) O_{p'}(J)= J\cap \la b\ra M \unlhd G.$$
Let $\la c\ra=P\cap M$. Clearly $P=\la a\ra \times \la b\ra \times
\la c\ra$. Let $$ C=\left \{\begin{array}{lrrrrrrlr} \la c\ra
\ltimes [O_{p'}(G), \la c\ra],
&{\rm if } \,\,G\,\,{\rm is} \,\,p-{\rm solvable},&\\
G^{{\rm Sol}_p}, &{\rm otherwise}.&
\end{array}\right.$$
Assume that $G$ is $p$-solvable. Since $G$ admits an abelian Sylow
$p$-subgroup with $G=O^{p'}(G)$, we have $G=P\ltimes O_{p'}(G)$.
Thus  $C=\la c\ra \ltimes [O_{p'}(G), \la c\ra]$ is normal in $G$,
and  $C=O^{p'}(C)$. Now let us investigate the case when  $G$ is
non-$p$-solvable. Clearly $C=G^{{\rm Sol}_p}$ is also normal in $G$,
and also $C=O^{p'}(C)$. Observe that $|C|_p=|G^{{\rm Sol}_p}|_p\geq
p$ because $G$ is not $p$-solvable, and that  $C=G^{{\rm Sol}_p}\leq
M$ because $G/M\cong E(p^2)$ is $p$-solvable. It follows  that
$|C|_p=|M|_p=p$ and $P\cap C=P\cap M=\la c\ra$. Now we always have
$$\la a\ra, B, C\unlhd G, \,\, P\cap B=\la b\ra,\,\, P\cap C=\la c\ra,$$
while $B=O^{p'}(B)$ and $C=O^{p'}(C)$. Now   $\la a\ra BC$ is normal
in $G$ with $p'$-index. Since $G=O^{p'}(G)$, we have
$$G=\la a\ra BC=\la a\ra \times (BC).$$
Since $J=\la a\ra \times B$ is solvable but $G$ is not solvable,
$C$ is necessarily nonsolvable.

Now it suffices to show that $B\cap C=O_{p'}(B)$.  Let $D=B\cap C$.
Note that $D$ is a $p'$-group because $|BC|_p=p^2=|B|_p|C|_p$.
Consequently $$|B/D|_p=|B|_p=p=|C|_p=|C/D|_p.$$ Assume that $C/D$
has a normal Sylow $p$-subgroup. Since $O^{p'}(C)=C$ and $|C/D|_p=p$,
we have  $|C/D|=p$. However $D\leq B< J\leq {\rm Sol}(G)$, we get
that  $C$ is solvable,  a contradiction. Consequently $C/D$ admits a
nonnormal Sylow $p$-subgroup. Assume that $B/D$ also admits a
nonnormal Sylow $p$-subgroup. By the Ito-Michiler Theorem, we may
take $\theta_b\in {\rm Irr}(B/D)$ and $\theta_c\in {\rm Irr}(C/D)$
such that $p$ divides $\theta_b(1)$ and $\theta_c(1)$. It follows
that $e_p(H)\geq e_p(G)\geq e_p(B/D\times C/D)\geq \log_p
((\theta_b\theta_c)(1))_p\geq 2$, a contradiction. Therefore, $B/D$
must admit a normal  Sylow $p$-subgroup. Since $B=O^{p'}(B)$ and
$|B/D|_p=p$, we have  $B/D\cong C_p$. Consequently  $D=O_{p'}(B)$,
as required.

\bigskip

(6) Final contradiction.

Since $J=O^{p'}(J)$ and $F(J)$ is cyclic by (4),
$J/F(J)=J/C_J(F(J))$ is an abelian  $p$-group. Since $|P_J|=p^2$ and
$|F(J)|_p=|O_p(G)|=p$, we have  $$J/F(J)\cong C_p.$$ We keep all
notations in (5). Let $D=O_{p'}(B)$ and let us investigate $BC\unlhd
G$. By (5), it  is easy to see that
$$B=\la b\ra \ltimes D,\,\,D=B\cap C=O_{p'}(B)=O_{p'}(J)=O_{p'}(F(J)).$$
Since $O_p(G)\cong C_p$, we have $D>1$ by (5b). Let $D/S$ be a
$G$-chief factor. Since $D=O_{p'}(F(J))$ is a cyclic $p'$-group,
$D/S\cong C_q$ for some prime $q\not=p$. Since $G=O^{p'}(G)$ and
$G/C_G(D/S)$ is abelian, all $p'$-elements of $G$ are contained in
$C_G(D/S)$, that is, $O^p(G)\leq C_G(D/S)$. Consequently
\begin{equation}\label{eq212}O^{p}(BC)\leq C_G(D/S)\cap BC= C_{BC}(D/S).
\end{equation}

Suppose that $D/S\leq \Phi(BC/S)$.
Since $B/S\unlhd BC/S$ and $\la b\ra S/S\in {\rm Syl}_p(B/S)$,
we conclude by   the Frattini Argument  that
$$BC/S=(B/S) (N_{BC/S}(\la b\ra S/S))
=(D/S) (N_{BC/S}(\la b\ra S/S))=N_{BC/S}(\la b\ra S/S).$$
It follows that $\la b\ra S/S\unlhd B/S$.
Consequently $B/S=\la b\ra S/S \times D/S$ and $O^{p'}(B)< B$, a contradiction.

Suppose that $D/S\cap \Phi(BC/S)=1$. Let $\lambda\in {\rm
Irr}^{\sharp}(D/S)$. We also view $\lambda$ as a character of $D$
and write  $T={\rm I}_{BC}(\lambda)$. Since   $BC/S$ is split over
$D/S$,  $\lambda$ is extendible to $T$. Since $D/S\cong C_q$, we
know that an element $g\in BC$ fixes $\lambda$ if and only if $g$
centralizes $D/S$. Consequently
\begin{equation}\label{eq213} T=C_{BC}(D/S).\end{equation}
Assume that $T\cap B\not=D$. Since $T\cap B\geq D$ and $B/D\cong
C_p$, we have $T\cap B=B$, that is, $B\leq T$. Now  $B$ centralizes
$D/S$ by (\ref{eq213}), and this also yields $O^{p'}(B)<B$, a
contradiction. Therefore $T\cap B=D$, and in particular $|T|_p\leq
p$. Note that $e_p(BC/S)\leq 1$, and it follows that $p^2$ does not
divide $\chi(1)$ for any irreducible constituent $\chi$ of
$\lambda^{BC}$. By  \cite[Theorem 6.11, Corollary 6.17]{I}, $T/D$
admits a normal abelian Sylow $p$-subgroup of order $p$. Observe
that
$$BC/O^{p}(BC)\leq E(p^2), \,\, T\cap B=D,$$ and
 $$
BC>T= C_{BC}(D/S)\geq O^p(BC)$$ by (\ref{eq212}) and (\ref{eq213}).
It follows that
$$T\unlhd BC,\, BT\unlhd BC$$
and that
$$  |BT|_p=\frac{|B|_p|T|_p}{|B\cap T|_p}=\frac{|B|_p|T|_p}{|D|_p}=p^2,$$
that is, $BT$ is a normal subgroup of $BC$ with $p'$-index.  Note
that $O^{p'}(BC)=BC$ because $G=\la a\ra \times BC$ and
$G=O^{p'}(G)$. Consequently, $BC=BT$ and
$$BC/D=(B/D)(T/D)=B/D\times T/D.$$ Since $O^{p'}(BC)=BC$, we also have
$O^{p'}(T/D)=T/D$. Recall that  $T/D$ admits a normal abelian Sylow
$p$-subgroup of order $p$, and it follows that $T/D\cong C_p$. Now
$BC/D\cong E(p^2)$. Since $D$ is cyclic, $BC$ and $G, H$ are
necessarily solvable, a contradiction.\pfend

\bigskip

\noindent \textbf{Proof of Theorem \ref{t101}.}
It follows directly from Theorem \ref{t206} and Lemma \ref{l201}(2).\pfend

\bigskip

\begin{rem}Let $H$ be a finite group with $|H/O_p(H)|_p=p^2$
and $e_p(H)=1$ for  an odd prime $p$,
assume that $|U/O_p(U)|_p\leq p$ for all proper quotient groups $U$ of $H$.
Let us investigate the structure of $H$.\end{rem}

Suppose that $p=3$ and ${\rm Alt}_7$ is involved in $H$.
By Lemma \ref{l201}(3), we have $O^{3'}(H)/O_3(H)$ $\cong {\rm Alt}_7$.
Suppose that $p\geq 5$ or,
$p=3$ and ${\rm Alt}_7$ is not involved in $H$.
Using the same arguments as in the proof of Theorem \ref{t206},
we obtain that

(1) $O_p(H)=\Phi(H)=1$, $|H|_p=p^2$, $|{\rm Sol}(H)|_p\geq p$;

(2) $H$ admits a unique minimal normal subgroup, say $V$; $V$ is an
elementary $r$-group for a prime $r\not=p$; $H=G\ltimes V$ for a
maximal subgroup $G$ of $H$; $O_p(G)\cong C_p$;

(3) Every $\lambda\in {\rm Irr}^{\sharp}(V)$  is fixed by a unique
subgroup of $G$ of order $p$.
It follows from Lemma \ref{l203}
that ${\rm Irr}(V)$ and $V$ are faithful  primitive $O^{p'}(G)$-module.\pfend

\bigskip

We do not find  an example of  non-$p$-solvable group $H$ with
$e_p(H)=1$ and $|H/O_p(H)|_p=p^2$ for a prime $p\geq 5$. Can one
construct such an example?

\bigskip
It would be also desirable for one to obtain a detailed structure of
finite solvable groups $G$ with $e_p(G)=1$ (or $e_p'(G)=1$) and
$|G/O_p(G)|_p=p^2$.

\section{On $p$-parts of Brauer character degrees}

We will prove Theorem \ref{t103} at the end of this section. We begin with listing some known results
about finite groups $G$ with $e'_p(G)\leq 1$.

\begin{lem}\label{l301}
Let $p$ be a prime and $G$ be a finite $p$-solvable group with $O_p(G)=1$.
Then $e'_p(G)\leq 1$ if and only if $e_p(G)\leq 1$.
\end{lem}
\pf This is \cite[Corollary 1.4]{LNTV}.

\bigskip

\noindent\textbf{ Proof of Corollary \ref{c102}.}
It follows directly from Theorem \ref{t101} and Lemma \ref{l301}.

\bigskip

\begin{lem}\label{l302}
Let $p$ be a prime,
let $G$ be a finite group with $O_p(G)=1$,
and assume that $G$ has
an abelian Sylow p-subgroup.
If $e'_p(G)\leq 1$, then $e_p(G)\leq 1$. \end{lem}
\pf This is  \cite[Theorem 1.3]{LNTV}.\pfend

\begin{lem}\label{l303}
Let $p$ be a  prime and let $G$ be a non-$p$-solvable group
with $G=O^{p'}(G)$  and ${\rm Sol}_p(G)=1$.
If  $e'_p(G)\leq 1$,
then $G$ is a nonabelian simple group,
furthermore the following hold:

{\rm (1)} If $p=2$, then $G\cong M_{22}$;

{\rm (2)} If $p\geq 3$, then either
$|G|_p=p$ or, $p=3$ and $G\cong {\rm Alt}_7$. \end{lem}
\pf  It follows directly from Lemmas 4.5 and 4.3 of \cite{LNTV}. \pfend

\bigskip

Let $p$ be a prime and $G$ be a finite group. Recall that $\chi\in
{\rm Irr}(G)$ is said to have $p$-defect zero if $\chi(1)_p=|G|_p$.
By \cite[Theorem 3.8]{N}, if $\chi$ is of $p$-defect zero, then
$\chi^0\in {\rm IBr}_p(G)$, where $\chi^0$ is the restriction of
$\chi$ to the set of $p$-regular elements of $G$.

\begin{pro}\label{p304} Let $G$ be a nonsolvable group.
Then $e'_2(G)\leq 1$
if and only if
$G/O_2(G)$ is a direct product of $M_{22}$ and an odd order group.\end{pro}

\pf Since ${\rm IBr}_2(G)={\rm IBr}_2(G/O_2(G)$, we may assume by
induction that $O_2(G)=1$. Since $e'_2(M_{22})=1$ by \cite{B-Atlas},
we  only need to prove the necessity. Note that all odd order groups
are solvable.

Suppose that $e'_2(G)\leq 1$ and let  $W=O^{2'}(G/{\rm Sol}(G))$.
Then $W$ is nonsolvable with $e'_2(W)\leq 1$. Clearly $W=O^{2'}(W)$
and ${\rm Sol}_2(W)={\rm Sol}(W)=1$. Hence $W\cong M_{22}$ by Lemma
\ref{l303}. Since ${\rm Out}(M_{22})\cong C_2$ by \cite{Atlas},
$G/{\rm Sol}(G)=W\times (K/{\rm Sol}(G))$, where $K/{\rm Sol}(G)$ is
of odd order. This implies that $K={\rm Sol}(G)$ and
\begin{equation}\label{eq301} G/{\rm Sol}(G)\cong M_{22}.\end{equation}
Assume that $G>G'$ and let $M\unlhd G$ be such that $|G:M|$ is a prime.
Note that $M$ is nonsolvable with  $e'_2(M)\leq 1$.
Since $O_2(M)\leq O_2(G)=1$,
we conclude by  induction that
$M=S\times U$ where $S\cong M_{22}$ and $U=O_{2'}(M)={\rm Sol}(M)$.
Note that if  $M=M {\rm Sol}(G)$, then
${\rm Sol}(G)\leq  M$, and this contradicts  (\ref{eq301}). Therefore
$$G=M {\rm Sol}(G)=SU{\rm Sol}(G)=S{\rm Sol}(G)=S\times {\rm Sol}(G).$$
Now $e'_2(G)\leq 1$ implies  $e'_2({\rm Sol}(G))=0$.
Therefore ${\rm Sol}(G)$ admits a normal Sylow $p$-subgroup.
Since  $O_2(G)=1$,
${\rm Sol}(G)$ has odd order,
and we are done. Hence we may assume that
\begin{equation}\label{eq302}
O_2(G)=1,\,G=G', G/{\rm Sol}(G)\cong M_{22}.
\end{equation}

Now it suffices to show that $G\cong M_{22}$.
Assume that this  is not true.
By (\ref{eq302}),
we may take $N\unlhd G$ maximal so that $G/N$ has an odd order chief factor.
By (\ref{eq302}) and the maximality of $N$,
we get that $N\leq {\rm Sol}(G)$,
$G/N$ has a minimal normal subgroup $V/N\cong
E(r^n)$ for an odd prime $r$,
 $V/N$ is the unique minimal normal subgroup of $G/N$,
and that ${\rm Sol}(G)/V=O_2(G/V)$.
In order to see a contradiction, we may assume that $$N=1.$$
Since $e'_2({\rm Sol}(G))\leq 1$,
we have $e_2({\rm Sol}(G))\leq 1$ by Lemma \ref{l301}.
This also implies by Lemma \ref{l201}(1) that
${\rm Sol}(G)/V$ is an abelian $2$-group of order at most $4$.
Note that  $e_2({\rm Sol}(G))=\log_2|{\rm Sol}(G)/V|$ by \cite[Lemma 18.1]{MW}.
It follows that $|{\rm Sol}(G)/V|\leq 2$.
Since $G=G'$, we conclude that  $$G/V\cong M_{22}\,\, {\rm or}\,\,2.M_{22}.$$

Assume that $G/V\cong M_{22}$. By \cite[Lemma 2.3]{LLQ17}, there
exists $\chi\in {\rm Irr}(G)$ of $2$-defect zero. This leads to
$\chi^0 \in {\rm IBr}_2(G)$. Now $e'_2(G)\geq \log_2
(\chi^0(1))_2=7$, and we get a contradiction.

Assume that $G/V\cong 2.M_{22}$.
Note that ``$V\leq \Phi(G)$'' would imply ${\rm Sol}(G)$ is nilpotent,
and this contradicts the unique minimal normality of $V$.
Consequently $\Phi(G)=1$,
whence $$G=H\ltimes V,\,\,{\rm where}\,\,H\cong G/V\cong 2. M_{22}.$$
Let $\lambda\in {\rm IBr}_2^{\sharp}(V)$
and  write $T={\rm I}_H(\lambda)$.
Clearly ${\rm IBr}_2(V)={\rm Irr}(V)$ because $V$ has odd order.
Note  that  ${\rm Sol}(G)=Z(H)\ltimes V$ is a Frobenius group with $V$ as its kernel.
It follows that
\begin{equation}\label{eq303}T\cap Z(H)=1.\end{equation}
Since $4$ does not divide $\phi(1)$ for any $\phi\in {\rm
IBr}_2(G|\lambda)$, we conclude by (\ref{eq303}) and Clifford
Correspondence (\cite[Theorem 8.9]{N})
 that
a Sylow $2$-subgroup $P$ of $T$ must be isomorphic
to a Sylow $2$-subgroup of $M_{22}$.
Now (\ref{eq303}) also implies that
$$P\times Z(H) \in {\rm Syl}_2(H).$$
Applying  Lemma \ref{l205},
we have  $Z(H)\cap H'=1$.
However $H=H'$ because $G=G'$, we get a contradiction.
\pfend

\begin{pro}\label{p305}
Let $G$ be a finite group and $p$ be an odd prime.
Assume that $e'_p(G)\leq 1$.
Then  $G/O_p(G)$ admits an elementary abelian
Sylow $p$-subgroup.\end{pro}
\pf By induction, we may assume that $O_p(G)=1$ and $G=O^{p'}(G)$.
Let $P\in {\rm Syl}_p(G)$
 and let $N={\rm Sol}_p(G)$.
Since $N$ is $p$-solvable and $e'_p(N)\leq 1$, we have $e_p(N)\leq
1$ by Lemma \ref{l301}. Hence $|P\cap N|\leq p^2$ by Theorem
\ref{t101}. Assume that $P\cap N$ is a cyclic group  of order $p^2$.
It follows by \cite{Z} that $N$ has an irreducible character of
degree divisible by $p^2$ , a contradiction. Hence $P\cap N\leq
E(p^2)$. Now we may assume $P>P\cap N$. Observe that $G/N$ is
non-$p$-solvable with $e'_p(G/N)\leq 1$, while $G/N=O^{p'}(G/N)$ and
${\rm Sol}_p(G/N)=1$. It follows from Lemma \ref{l303} that $G/N$ is
a nonabelian simple group, and that either $|G/N|_p=p$, or $p=3$ and
$G/N\cong {\rm Alt}_7$. In particular, $G/N$ always has an
elementary abelian Sylow $p$-subgroup. Now  the required result
holds obviously for the case when $P\cap N=1$. Consequently, we may
assume that
$$P> P\cap N \cong C_p \,\, {\rm or}\,\,E(p^2).$$

Let $V$ be a minimal normal subgroup of $G$. Assume that $V\nleq N$.
Since $G/N$ is a nonabeian simple group, we have  $G=N\times V$
where $V\cong G/N$. Since both $V$ and $N$ have elementary  abelian
Sylow $p$-subgroups, $P$ is elementary abelian, and we are done. We
therefore may assume that $V\leq N={\rm Sol}_p(G)$. Now $V\leq
O_{p'}(G)$ because $O_p(G)=1$. Note that if $O_p(G/V)=1$, then the 
induction implies the required result. Hence we may also assume that
$D/V:=O_p(G/V)>1$.

Let us investigate the normal subgroup $D$ of $G$. Clearly  $D\leq
N$ and $O_p(D)=1$. Since $P\cap N\leq E(p^2)$, we have  $C_p\leq
D/V\leq E(p^2)$. By \cite[Corollary 2.8]{Q12} when $V$ is nonabelian
or \cite[Lemma 18.1]{MW} when $V$ is abelian, we may take
$\lambda\in {\rm Irr}^{\sharp}(V)={\rm IBr}_p^{\sharp}(V)$ such that
${\rm I}_D(\lambda)=V$, that is, ${\rm I}_G(\lambda)\cap D=V$. Since
$p^2$ does not divide $\phi(1)$ for any $\phi\in {\rm
IBr}_p(G|\lambda)$, we conclude by \cite[Theorem 8.9]{N} that
$$D/V\cong C_p$$ and that ${\rm I}_G(\lambda)$
contains a Sylow $p$-subgroup $Y$ of order $|P|/p$. Note that $D/V$
is central in every Sylow $p$-subgroup of $G/V$. This implies that
$$(YV/V)\times (D/V)\in {\rm Syl}_p(G/V).$$ Applying Lemma
\ref{l205}, we get that $(D/V)\cap (G/V)'=1$, that is,
\begin{equation}\label{eq304} D\cap G'V=V.\end{equation}
Since  $G/N$ is a nonabelian simple group,
we have  $G=NG'$.
Observe  that $G/G'$ is a $p$-group because $G=O^{p'}(G)$
and that $P\cap N$ is elementary abelian.
It follows that
$G/G'\cong N/(N\cap G')$ is an elementary abelian $p$-group.
Now all subgroups of $G/G'$ are complemented in $G/G'$.
Let $R/G'$ be a complement of $DG'/G'$ in $G/G'$.
We have
\begin{equation}\label{eq305}
G =RD,\,\, R\cap DG'=G'.\end{equation}
Using (\ref{eq304}), (\ref{eq305}) and the Dedekind's Modular Law,
we get that
$$D\cap RV=V(D\cap R)=V(D\cap DG'\cap R)
=V(D\cap G')=V(D\cap G'V\cap  G')=V(V\cap G')=V,$$
and that
\begin{equation}\label{eq306}
G/V=RD/V=(D/V) \times (RV/V).
\end{equation}
Now $RV$ is a proper normal subgroup of $G$ with  $O_p(RV)\leq
O_p(G)=1$. Applying the inductive hypothesis to $RV$, we conclude
that $RV$ has an elementary abelian Sylow $p$-group. This implies by
(\ref{eq306}) that $G/V$ has an elementary abelian Sylow
$p$-subgroup. Therefore $P\cong PV/V$ is elementary abelian. \pfend

\bigskip

Suppose that $e_p(G)=1$ for an odd prime $p$. By Theorem \ref{t101}
and a result in \cite{Z},   we also conclude that $G/O_p(G)$ admits
an elementary abelian Sylow $p$-subgroup of order at most $p^2$.

\begin{cor}\label{c206}
Let $G$ be a finite group with $O_p(G)=1$ for an odd prime $p$.
If $e_p'(G)\leq 1$, then $e_p(G)\leq 1$. \end{cor}

\pf It follows from Proposition \ref{p305} and Lemma
\ref{l302}.\pfend

\bigskip

\noindent{\textbf{Proof of Theorem \ref{t103}.} If $p=2$, the
required result follows from Corollary \ref{c102} and Proposition
\ref{p304}. For the $p\geq 3$ case, we have $e_p(G/O_p(G))\leq 1$ by
Corollary \ref{c206}. Now Theorem \ref{t101} implies the
result.\pfend

\bigskip

\begin{center} {\bf ACKNOWLEDGMENT} \end{center}

The authors would like to thank Dr. Yanjun Liu for his helpful discussions.

\end{document}